\documentclass{IEEEtran}

\usepackage{balance} 

\usepackage{xcolor,calc}

\usepackage[pdftex]{graphicx}
\usepackage{amsmath} 
\usepackage{amssymb}  
\usepackage[noadjust]{cite}
\usepackage{color}
\usepackage{enumerate}
\usepackage{subfigure}
\usepackage{multirow}
\usepackage{mathtools}
\usepackage{booktabs}
\usepackage{epstopdf}

\usepackage[normalem]{ulem}

\usepackage{hyperref}
\usepackage{url}
\usepackage{bm}

\usepackage{booktabs}
\newcommand{\ra}[1]{\renewcommand{\arraystretch}{#1}}

\hypersetup{colorlinks,citecolor=black,filecolor=black,linkcolor=black,urlcolor=black}
\newtheorem{theorem}{Theorem}

\newtheorem{definition}{Definition}
\newtheorem{assumption}{Assumption}

\usepackage{chngcntr}

\let\origtheassumption\theassumption

\title{\LARGE \bf
A Non-Cooperative Game Approach to Autonomous Racing
}

\author{Alexander Liniger and John Lygeros
\thanks{The authors are with the Automatic Control Laboratory, ETH Z\"urich, 8092 Z\"urich, Switzerland (e-mails: 
liniger@control.ee.ethz.ch; lygeros@control.ee.ethz.ch)

This manuscript is the preprint of a paper submitted to IEEE Transaction on Control System Technology and is subject to IEEE copyright. IEEE maintains the sole rights of distribution or publication of the work in all forms and media. If accepted, the copy of record will be available at {http://ieeexplore.ieee.org/}}
}

\begin{document}

\maketitle
\thispagestyle{empty}
\pagestyle{empty}

\begin{abstract}
We consider autonomous racing of two cars and present an approach to formulate racing decisions as a non-cooperative non-zero-sum game. We design three different games where the players aim to fulfill static track constraints as well as avoid collision with each other; the latter constraint depends on the combined actions of the two players. The difference between the games are the collision constraints and the payoff. In the first game collision avoidance is only considered by the follower, and each player maximizes their own progress towards the finish line. We show that, thanks to the sequential structure of this game, equilibria can be computed through an efficient sequential maximization approach. Further, we show these actions, if feasible, are also a Stackelberg and Nash equilibrium in pure strategies of our second game where both players consider the collision constraints. The payoff of our third game is designed to promote blocking, by additionally rewarding the cars for staying ahead at the end of the horizon. We show that this changes the Stackelberg equilibrium, but has a minor influence on the Nash equilibria. For online implementation, we propose to play the games in a moving horizon fashion, and discuss two methods for guaranteeing feasibility of the resulting coupled repeated games. Finally, we study the performance of the proposed approaches in simulation for a set-up that replicates the miniature race car tested at the Automatic Control Laboratory of ETH Z\"urich. The simulation study shows that the presented games can successfully model different racing behaviors and generate interesting racing situations.
\end{abstract}

\section{Introduction}
One of the major challenges in self-driving cars is the interaction with other cars. This holds true for public roads as well as race tracks, as the fundamental problem comes from the fact that the intentions of the other cars are not known. Autonomous driving and autonomous racing have been successfully demonstrated repeatedly \cite{Horowitz2000,Buehler2007,stanford,Liniger_2014}, while driving in close proximity with unpredictable vehicles on public roads was addressed in the DARPA Urban Challenge \cite{Buehler2009}. Urban driving remained an active field of research ever since \cite{Ziegler2014,Paden2016}. Here we stay in this realm of driving in the presence of non-cooperative neighbors but concentrate on the more controlled environment of autonomous racing on a race track, where the rival cars/opponents generate a dynamically changing environment. In autonomous racing the goal is to drive as fast as possible around a predefined track, making it necessary to control the car at the limit of handling \cite{stanford}. Moreover, in racing the interactions between the cars are governed by less strict rules than on public roads; for example, there are no lanes and lane changes are not indicated. This implies that while in autonomous driving one can develop methods that assume a minimum degree of cooperation from neighboring vehicles \cite{Althoff2009,Carvalhoa2014,Carvalhoa2015}, such an approach is more difficult for autonomous racing. Therefore, methods for autonomous racing often ignore obstacle avoidance altogether \cite{stanford,Verschueren14,Rosolia2016}, or only consider the avoidance of static obstacles \cite{Liniger_2014,Funke2017}.

To deal with the uncertainty in the driving behavior of the opponent cars, we propose here a game theoretical approach to autonomous racing. Zero-sum game methods have been investigated for related problems in air traffic \cite{Tomlin1998} and autonomous driving \cite{Lygeros1998}. Autonomous racing could be formulated as a reach-avoid game \cite{Margellos2011}, where efficient computational methods for open-loop games exit \cite{Zhou_2012,Zhou_2018}. However, since such games assume the worst-case behavior of the opponent, special care is needed to ensure that the resulting driving style is not too defensive for autonomous racing. Game theory has also been used to derive driver models to simulate and verify autonomous driving algorithms \cite{Li2016,Yoo2012}. Those approaches do not assume a worst-case behavior and are typically investigated by means of Stackelberg equilibria. 

Adopting some of the assumptions of these game theory based driver models, we assume that in a race the competing car has no benefit from causing a collision. We propose to model the interactions between race cars as a non-cooperative non-zero-sum game, where the players only get rewarded if they do not cause a collision. By restricting our attention to finite horizon two-player games, we derive three different games that, under the additional hypothesis that the action set of the two players is finite are formulated as bimatrix games. In the first game, both players maximize their progress, but only the follower is concerned with avoiding collisions. For this game, we show that the Stackelberg and Nash equilibria can be computed by a sequential maximization approach, where the leader determines his trajectory independently of the follower, and the follower plays his best response. In the second game, both players consider collisions. This allows to prove that the players always choose a collision-free trajectory pair if possible. Further, we show that if the sequential maximization approach leads to a collision-free trajectory pair, this trajectory pair is also a Stackelberg and Nash equilibrium of the second game. In the third game, we propose a modified payoff that promotes blocking behavior, by rewarding actions that do not necessarily maximize progress towards the finishing line but instead aim at preventing the opponent from overtaking. For the third game, we show that the Stackelberg equilibrium is a blocking trajectory pair if one exists, which stands in contrast to the Nash equilibrium which only results in blocking trajectories in particular cases.

The approach proposed here is related to \cite{Williams2017}, where a similar game theoretic approach to autonomous racing is proposed. Also related is the approach proposed in \cite{Spica2018}, which formulates a zero-sum game for drone racing, which is played in a receding horizon fashion. However, in both \cite{Williams2017} and \cite{Spica2018} the equilibria of the games are not analyzed and only best-response dynamics are consider to solve the game. One could in principle also consider extending autonomous racing approaches with static obstacle avoidance, such as those developed in \cite{Liniger_2014,Funke2017}, to incorporate dynamical obstacles by predicting the opponent movement and using time-varying constraints. This, however, would not allow the flexibility of our game theoretic approach, where both players influence each other's decisions.

For online implementation, we propose repeating the game in a receding horizon fashion, similar to Model Predictive Control (MPC), giving rise to a sequence of coupled games. The use of moving horizon games has been proposed in other applications \cite{Hespanha2000,Cruz2002,Virtanen2006,Takei_2012}, as a method to tackle higher dimensional problems not tractable by dynamic programming. We investigate the sequence of coupled games and propose modified constraints based on viability theory which guarantee recursive feasibility and introduce an exact soft constraint reformulation that guarantees feasibility at all times. Compared to \cite{Liniger_2017} where viability theory was used to guarantee recursive feasibility only with respect to track constraints, here recursive feasibility is additionally guaranteed with respect to the actions of the opponent player.

Finally, we investigate the game formulations in a simulation study, replicating the miniature race car set-up presented in \cite{Liniger_2014}, and investigate the influence of the game formulation on different performance indicators such as the number of overtaking maneuvers and the collision probability.

This paper is structured as follows, in Section \ref{sec:GameIngre} we introduce the ingredients to formulate the racing games. Three different racing games are formulated in Section \ref{sec:GameForm}. The optimal solution in terms of the Stackelberg and Nash equilibrium is studied in Section \ref{sec:Equilibrium}. In Section \ref{sec:MovingHorizonGames}, we discuss feasibility of moving horizon games. The performance of the approaches is investigated in a simulation study in Section \ref{sec:SimResults}, while Section \ref{sec:Conclusion} provides some concluding remarks. Finally, in Appendix \ref{app:ppModel} we give a summary of the vehicle model used in this paper, as presented in \cite{Liniger_2017}.

\section{Game ingredients}\label{sec:GameIngre}

In a car race, the goal is to finish first; mathematically this can be interpreted as a non-cooperative dynamic game where each player tries to reach the finishing line before any other player. In \emph{Formula One} or similar racing series this game is tackled by separating it into two problems: the first problem concerns slowly changing strategic decisions (pit stops, power consumption, tire wear, etc.) and is solved by the pit crew \cite{Bekker2009}. These high-level decisions are then executed by a skilled race car driver, who solves the second problem by driving the car at the handling limit while interacting with other cars. In this paper, we concentrate on the task of the race car driver and formulate it as a dynamic game. To formulate the racing games three ingredients are required, first an appropriate vehicle model, second an objective function representing/approximating the ``finish first" goal and third state constraints characterizing collisions.

\subsection{Dynamics}\label{sec:GameIngreDyn}
We start with the path planning model of \cite{Liniger_2017} as the vehicle model, where the dynamics are described by a nonlinear discrete-time system of the form
\begin{align*}
\bar{x}_{k+1}=\bar{f}(\bar{x}_k,u_k).
\end{align*}
The state $\bar{x}=[X, Y, \varphi, q] \in \mathbb{R}^4$ represents the $X$ and $Y$ position in global coordinates, the heading angle $\varphi$ and $q$ denotes the constant velocity mode or in other words the current motion primitive. The input $u\in U(q)$ decides which constant velocity mode is chosen next, where $U(q)$ represents concatenation constraints on the constant velocity modes, reflecting for example acceleration limits. The update equation for the state $\bar{f}: \mathbb{R}^4 \times \mathbb{R} \rightarrow \mathbb{R}^4$ encodes the kinematic constraints of the bicycle model and the dynamic constraints of the constant velocity segment updates and is defined in detail in \cite{Liniger_2017}. We note that strictly speaking the model is hybrid since one of the states, the index $q$ of the constant velocity segment, takes values in a finite set. Since the input constraints $U(q)$ ensure that if $q$ is initialized in this finite set it will never leave it, to simplify notation we can ignore this complication, however, embed $q$ in the real line, and think of it as real-valued. For a more detailed description of the model see Appendix \ref{app:ppModel}.

\subsection{Payoff}
To define the objective function of the game, we extend the state by an additional variable, a lap counter $c$, leading to $x=[\bar{x}, c] \in \mathbb{R}^5$.  To define the dynamics of the lap counter we consider the center line of the track as a parametrized curve $\mbox{\em cl}:[0, L) \rightarrow \mathbb{R}^2$, where $0$ represents the ``starting line"; note that the argument of the function  $\mbox{\em cl}$ loops from $L$ back to $0$ every time the car completes a lap. Given a state $x \in \mathbb{R}^5$, we define its projection $\bar{p}:\mathbb{R}^5\rightarrow [0,L)$ onto the center line using the first two components of the state
\begin{align} \label{eq:proj}
\bar{p}(x)=\arg\min_{l\in [0, L)}\|\mbox{\em cl}(l)-(X,Y)\|_2\,.
\end{align}
Under the reasonable assumption that the track length is such that the car can only cover a fraction of one lap in a single sampling time, we then define the update equation of the lap counter by
\begin{align*}
c_{k+1}=\left\{\begin{array}{ll}
c_k+1 & \mbox{ if } \bar{p}(x_k) > \bar{p}([\bar{f}(\bar{x}_k, u_k),c_k]) \\
c_k & \mbox{ otherwise}\,,
\end{array}\right.
\end{align*}
initialized with $c_0 = 0$. Appending this update equation to the dynamics $\bar{f}$ of $\bar{x}$ gives rise to the dynamics of the combined state
\begin{align*}
x_{k+1}=f(x_k, u_k)\,.
\end{align*}
Note that the lap counter also takes values in the discrete set of integers, but since its update equation is such that if initialized in this set it remains in this set, we again embed in the real line to simplify the notation. Finally, given a state $x$ we define the progress function
\begin{align} \label{eq:prgress}
p(x)=\bar{p}(x)+cL\,, 
\end{align}
that will form the basis for the objective function defined below.

To encode the objective function of each car as well as the fact that the cars should avoid collisions with each other one needs to consider both cars; we use the superscript $\varrho=1, 2$ to distinguish the state of the two cars. Furthermore, to formulate the racing games both cars start at their initial state denoted by $x^{\,1}$ and $x^{2}$, and even thought both cars use the path planning model they may be described by different constant velocity modes, highlighted by the superscript in $f^{\varrho}(\cdot,\cdot)$ and $U^{\varrho}$. Therefore the dynamics of the two cars are given by,
\begin{align} 
&x^{\varrho}_0 = x^{\varrho}\,, \quad \varrho = 1,2\,, \quad k = 0,...,N-1\,,\nonumber\\
&x^{\varrho}_{k+1} = f^{\varrho}(x^{\varrho}_k,u^{\varrho}_k) \,, \quad u_k^{\varrho} \in U^{\varrho}(q_k^{\varrho}) \,. \label{eq:dynCon}
\end{align}
Due to the discrete nature of the admissible inputs $U^{\varrho}(q^{\varrho})$ there exist a finite number of state-input trajectories. These different trajectories are denoted by a subscript $i = 1,..., n$ for car 1 and $j = 1,..., m$ for car 2 and the time step along the trajectory with $k=1,...,N$, resulting in the notation ${x}^{1}_{i}(k)$ and ${x}^{2}_{j}(k)$ denoting the different trajectories of the two cars at time step $k$. 

Each trajectory of the two cars has an associated progress payoff, which is the progress \eqref{eq:prgress} of its final state. The progress payoff of the $i$-th and $j$-th trajectory of car 1 respectively car 2 is therefore given by $p(x^1_i(N))$ and $p(x^2_j(N))$. In addition to progress, the objective function of the game should also reflect the fact that the car would like to remain on the track and avoid collisions with other cars. The first requirement can be encoded by a set $\mathcal{X}_{\text{T}} \subseteq \mathbb{R}^5$ that encodes the requirement that the first two components of the state are within the physical limits of the track. One can then impose a constraint that for all $k = 1,...,N$, $x_i^1(k) \in \mathcal{X}_{\text{T}}$, and $x_j^2(k) \in \mathcal{X}_{\text{T}}$ (or more precisely, penalize the car in the objective function whenever this constraint is violated).

Finally the collision constraints couple the two cars, in the sense that they depend on the actions of both players. One can encode this constraint through the signed distance function\footnote{The singed distance between two sets $R^1 \in \mathbb{R}^2$ and $R^2\in \mathbb{R}^2$ is defined as $sd(R^1,R^2) = \text{dist}(R^1,R^2) - \text{pert}(R^1,R^2)$, with $\text{dist}(R^1,R^2) = \inf \{\| T \| | (T+R^1) \cup R^2 \neq \emptyset\}$ and $\text{pert}(R^1,R^2) = \inf\{\| T \| | (T+R^1) \cup R^2 = \emptyset\}$.} \cite{Schulman2013}. If we assume that each car can be described by a rotated rectangle (indeed any fixed shape) centred at the $(X, Y, \varphi)$ component of the state, then the signed distance of the corresponding sets induces a function $d: \mathbb{R}^5\times \mathbb{R}^5\rightarrow \mathbb{R}$ such that $d(x^1, x^2) \geq 0$ if and only if the two cars do not physically overlap. Thus, a trajectory pair $(i,j)$ does not have a collision if 
\begin{align} 
d(x_i^1(k),x_j^2(k)) \geq 0\,, \quad k = 1,...,N \,. \label{eq:colCon}
\end{align} 

\section{Game formulation}\label{sec:GameForm}
In this section we combine the game ingredients discussed in Section \ref{sec:GameIngre} to formulate bimatrix games representing the racing objective. Let us assume that the two cars are the players; player 1 (P1) has $n$ possible trajectories, and player 2 (P2) has $m$ possible trajectories to choose from, constructed by combinatorially selecting inputs from the corresponding sets $U^1(q^1)$ and $U^2(q^2)$. The bimatrix game has two $m \times n$ payoff matrices $A$ and $B$ with elements $a_{i,j}$ and $b_{i,j}$ representing, respectively, the payoffs of P1 and P2, if P1 adopts trajectory $i$ and P2 trajectory $j$. The payoff associated with a trajectory pair $(i,j)$ encodes information about both the progress and the constraints.

Before formally introducing our different racing games, let us give two definitions and an assumption used in the rest of the paper. 
\begin{definition} \label{def:feas}
A trajectory pair $(i,j)$ is \emph{feasible} if $x_i^1(k) \in \mathcal{X}_{\text{T}}$, $x_j^2(k) \in \mathcal{X}_{\text{T}}$, and $d(x_i^1(k),x_j^2(k)) \geq 0$ for all $k = 1,...,N$.
\end{definition}
\begin{definition}
The leader of the game is the car which is ahead at the beginning of the game, in other words P1 if $p(x^1(0)) \geq p(x^2(0))$ and P2 otherwise.
\end{definition}
Without loss of generality we assume that P1 is the leader of the game.
\begin{assumption} \label{ass:leader}
At the beginning of the game P1 is ahead of P2 or in other words $p(x^1(0)) \geq p(x^2(0))$.
\end{assumption}

In following we discuss three different approaches to combine the progress payoff, collision, and track constraints to formulate three different bimatrix games, with increasingly complex interactions between the players. In the first two games, the winning objective is encoded as ``drive as fast as possible while avoiding collision with other cars". In the third game, a more complex payoff structure is used to promote blocking behavior, or in other words prevent an overtaking maneuver by driving sub-optimally with respect to progress along the track. For simplicity, we consider racing between only two cars. Racing with multiple opponents can be treated similarly but is considerably more challenging from a computational point of view.

\subsection{Sequential Game}
In the first game, the payoff matrices are designed such that the players stay on track, drive as fast as possible and avoid collisions. Roughly speaking, feasible trajectory pairs receive a payoff equal to their progress  if a trajectory leaves the track, this trajectory receives a constant payoff of $\kappa$ strictly smaller than the payoff of any trajectory which remains within the track. Likewise, trajectory pairs that collide also receive a low payoff. However, due to the leader-follower structure present in racing, we assume that only the follower (P2) is concerned with collision avoidance. Therefore, if a trajectory pair has a collision, P2 receives a low payoff of $\lambda$, whereas P1 receives the payoff as though the collision was not present. The resulting payoff matrices can be computed as follows,

\begin{align} \label{eq:SeqRacingGame}
a_{i,j} &= \begin{cases}
\kappa &\;\; \text{if } \exists k\in \{1,...,N \}  : x_i^1(k) \not\in \mathcal{X}_{\text{T}}\\
p(x_i^1(N)) &\;\;\text{else}\\
\end{cases}\,, \nonumber\\[-0.1cm]
\\
b_{i,j} &= \begin{cases}
\kappa & \;\; \text{if } \exists k \in \{1,...,N \}  : x_j^2(k) \not\in \mathcal{X}_{\text{T}}\\
\lambda & \begin{array}{ll}
\text{else if } & \exists k\in \{1,...,N \}   : \\
& d(x_i^1(k),x_j^2(k)) < 0
\end{array}\\
p(x_j^2(N)) & \;\;\text{else} \\
\end{cases}\nonumber\,.
\end{align}
Intuitively, if $\kappa$ and $\lambda$ are smaller than $\min ( \min_{i}  ( p(x_i^1(N)),\min_{j}   p(x_j^2(N)))$ both player, if possible, choose a trajectory which remains within the track and P2 will avoid a collision, see Section \ref{sec:Equilibrium} for a formal discussion. This condition can be fulfilled by choosing $\kappa, \lambda < 0$, as the progress is always positive. For the rest of the paper, we assume that $0 > \lambda \geq \kappa$, which implicitly penalizes collisions less than leaving the track. 

Note that in this game the payoff of the leader does not depend on the actions of the follower; in other words, all elements of a row of the payoff matrix $A$ will have the same value. We denote by $a_i$ the leader payoff for action $i\in \Gamma_1$ (in other words, the value of all elements of row $i$ in matrix $A$) to highlight the fact that it does not depend on the actions of the follower. Due to this purely sequential structure the game is called the \emph{sequential game}.

\begin{figure}[ht]
\centering
\includegraphics[width = 0.475\textwidth]{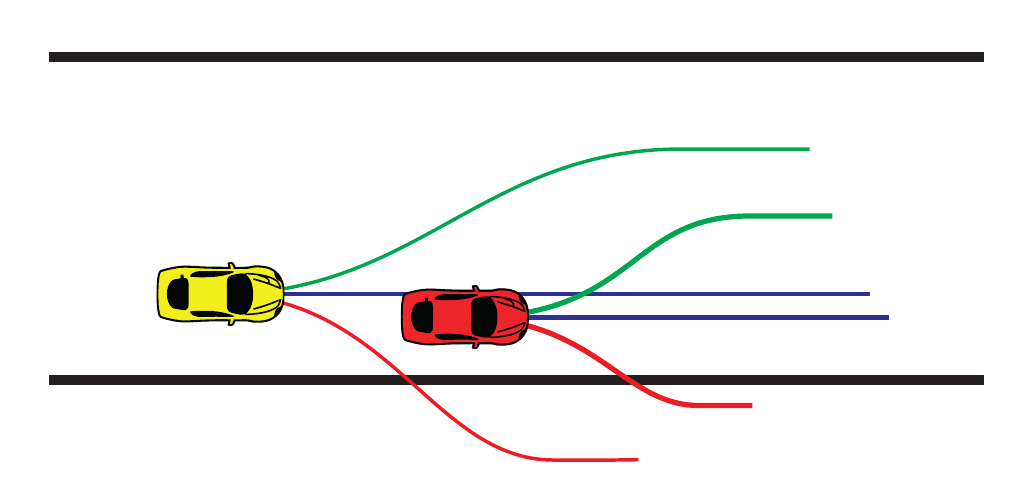}
    \caption{The leader P1 (red car) is in front of P2 (yellow car). Each player has three trajectories green, blue, red from top to bottom, and the trajectories are numbered from top (1) to bottom (3).}\label{fig:raceSit}
\end{figure}

To visualize the sequential game \eqref{eq:SeqRacingGame} let us look at one racing situation with only three possible trajectories for each player, see Fig. \ref{fig:raceSit}. We determine the progress by assuming that the zero point is at the beginning of the track interval shown and that the interval is 0.95\,m long, with 12\,cm long cars, which corresponds to the test track in our lab \cite{Liniger_2014} and the scale used in the subsequent simulation study. We use $\kappa = -10$ and $\lambda = -1$ leading to the payoff matrices $A$ (for the leader P1) and $B$ (for P2) of the racing game \eqref{eq:SeqRacingGame}

\begin{align} \label{eq:ABraceSitSeq}
A = \begin{bmatrix}
0.83 & {0.83} &0.83\\
{0.88} & 0.88 & 0.88 \\
-10 & -10 &  -10
\end{bmatrix} \, B = \begin{bmatrix}
0.81 & {0.86} &-10\\
{0.81} & -1 & -10 \\
0.81 & 0.86 &  -10
\end{bmatrix}.
\end{align}

For each player, the bottom trajectory leaves the track, so the payoff is $\kappa = -10$ for the corresponding trajectory, irrespective of the actions of the other player. This leads to the third row of $A$ and third column of $B$ identically equal to $-10$. If both players select their middle trajectory, there is a collision, giving rise to the payoff $\lambda =-1$ in the $(2,2)$ entry of $B$, note that $a_{2,2}$ is equal to the progress. The remaining combinations of actions neither leave the track nor cause a collision. Hence the payoff is the progress, giving rise to the entries in the order of $0.8$ in the two matrices.

\subsection{Cooperative Game}
Compared to the sequential game in the cooperative game, the leader does also consider collisions. This is incorporated into the game by modifying the payoff such that also the leader receives a payoff of $\lambda$ if there is a collision. In other words, both players only receive the progress payoff if a trajectory pair is feasible. The resulting payoff matrices can be computed as follows,
\begin{align} \label{eq:racingGame}
a_{i,j} &= \begin{cases}
\kappa &\;\; \text{if } \exists k\in \{1,...,N \}  : x_i^1(k) \not\in \mathcal{X}_{\text{T}}\\
\lambda & \begin{array}{ll}
\text{else if } & \exists k\in \{1,...,N \}   : \\
& d(x_i^1(k),x_j^2(k)) < 0
\end{array}\\
p(x_i^1(N)) &\;\;\text{else}\\
\end{cases}\,, \nonumber\\[-0.1cm]
\\
b_{i,j} &= \begin{cases}
\kappa & \;\; \text{if } \exists k \in \{1,...,N \}  : x_j^2(k) \not\in \mathcal{X}_{\text{T}}\\
\lambda & \begin{array}{ll}
\text{else if } & \exists k\in \{1,...,N \}   : \\
& d(x_i^1(k),x_j^2(k)) < 0
\end{array}\\
p(x_j^2(N)) & \;\;\text{else} \\
\end{cases}\nonumber\,,
\end{align}

Note that if we choose $\kappa, \lambda < 0$ as suggested in the sequential game, in case of the cooperative game the optimal trajectory pair will be feasible whenever possible since none of the players would benefit from violating the constraints, see Theorem \ref{theo:Link} for a formal discussion. Due to this cooperative obstacle avoidance approach, we call it the \emph{cooperative game}. 

When investigating the racing situation in Fig. \ref{fig:raceSit}, the payoff matrices $A$ and $B$ for the cooperative game are
\begin{align} \label{eq:ABraceSit}
A = \begin{bmatrix}
0.83 & {0.83} &0.83\\
{0.88} & -1 & 0.88 \\
-10 & -10 &  -10
\end{bmatrix} \, B = \begin{bmatrix}
0.81 & {0.86} &-10\\
{0.81} & -1 & -10 \\
0.81 & 0.86 &  -10
\end{bmatrix}.
\end{align}
Compared to the payoff matrices of the sequential game \eqref{eq:ABraceSitSeq} the only difference is $a_{2,2}$ which is now $-1$, since P1 also considers collisions. 

\subsection{Blocking game} \label{sec:Blocking}
When racing, maximizing progress is only a means to an end, the real objective is to finish first. In some cases, it may, therefore, be beneficial to prevent the opponent from overtaking, even if that means less progress. We refer to this behavior as ``blocking". In our third game, we modify the payoff matrices of the cooperative game \eqref{eq:racingGame} to reward staying ahead at the end of the horizon, giving rise to the so-called \emph{blocking game}. More precisely we add a parameter $w \geq 0$, which allows to trade-off staying ahead and maximizing progress. The resulting payoff matrices can be computed as follows,
\begin{align} \label{eq:blockingGame}
a_{i,j} &= \begin{cases}
\kappa &\;\; \text{if } \exists k\in \{1,...,N \}  : x_i^1(k) \not\in \mathcal{X}_{\text{T}}\\
\lambda & \begin{array}{ll}
\text{else if } & \exists k\in \{1,...,N \}   : \\
& d(x_i^1(k),x_j^2(k)) < 0
\end{array}\\
p(x_i^1(N)) &\;\;\text{else if   } \;\;\; p(x_i^1(N)) < p(x_j^2(N))\\
p(x_i^1(N)) + w &\;\;\text{else}
\end{cases}\,, \nonumber\\[-0.1cm]
\\
b_{i,j} &= \begin{cases}
\kappa & \;\; \text{if } \exists k \in \{1,...,N \}  : x_j^2(k) \not\in \mathcal{X}_{\text{T}}\\
\lambda & \begin{array}{ll}
\text{else if } & \exists k\in \{1,...,N \}   : \\
& d(x_i^1(k),x_j^2(k)) < 0
\end{array}\\
p(x_j^2(N)) &\;\;\text{else if } \;\;\; p(x_i^1(N)) \geq p(x_j^2(N))\\
p(x_j^2(N)) + w & \;\;\text{else} \\
\end{cases}\nonumber\,.
\end{align}

Note that in this case there are two terms linking the payoff of the two players, the collision constraint in the second line in \eqref{eq:blockingGame} and the blocking reward in the fourth line.

\begin{figure}[h]
\centering
\includegraphics[width = 0.475\textwidth]{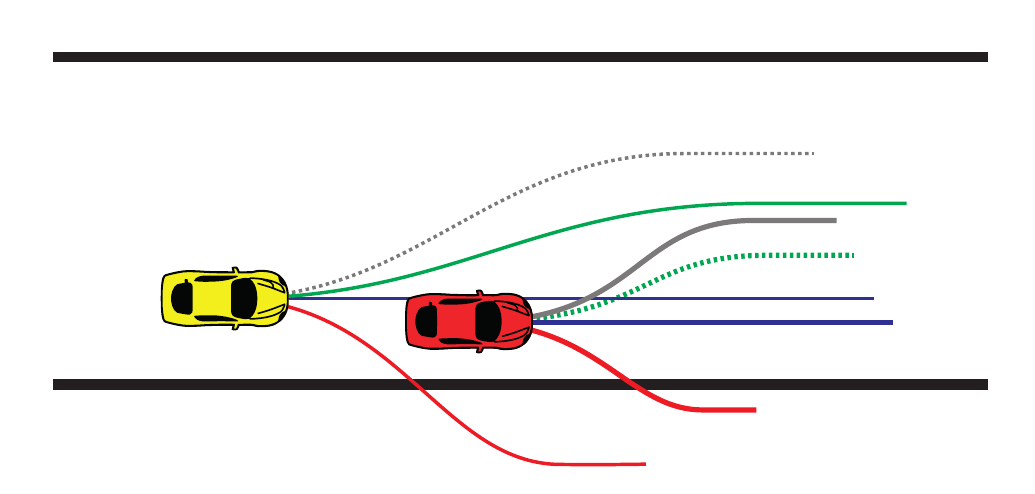}
    \caption{Illustration of a blocking situation, where P2 (yellow car) has a chance to overtake P1 (red car) at the end of the horizon by choosing the green trajectory (2). However, by choosing the dashed green blocking  trajectory P1 can avoid an overtaking maneuver, forcing P2 to choose his gray dashed trajectory (1) to avoid a collision. The trajectories are again numbered from top to bottom.}\label{fig:blocking}
\end{figure}

To visualize the effect of the additional payoff term, we look at a slightly more complex racing situation illustrated in Fig. \ref{fig:blocking}, where compared to the first racing situation in Fig. \ref{fig:raceSit} both players have one additional trajectory and P2 has the chance to overtake P1.

The payoff matrices of the bimatrix game corresponding to the racing situation in Fig. \ref{fig:blocking}, with $w = 0.5$ are,
\begin{align} \label{eq:BlockingEx}
&A = \begin{bmatrix}
0.83 + 0.5 &  -1  &0.83&0.83+0.5\\
0.85+0.5 & -1 &-1&0.85+0.5\\
{0.88+0.5} & 0.88 &-1 & 0.88+0.5\\
-10 & -10 &  -10&-10
\end{bmatrix} \,, \nonumber\\
&B = \begin{bmatrix}
0.81& -1 & {0.86} + 0.5 &-10\\
0.81& -1 & -1 &-10\\
{0.81}& 0.9 +0.5 & -1 & -10 \\
0.81+0.5& 0.9+0.5 & 0.86+0.5 &  -10
\end{bmatrix}.
\end{align}

\section{Characterization of equilibria} \label{sec:Equilibrium}

\subsection{Equilibrium concepts and feasibility}
Given the bimatrix games \eqref{eq:SeqRacingGame}, \eqref{eq:racingGame}, and \eqref{eq:blockingGame}, the question is how to find the optimal trajectory pair. Note that we only consider pure strategies, where the action space of the two players is given by  $\Gamma^1 := \{1,2,..,n\}$ and $\Gamma^2 := \{1,2,..,m\}$. Further, we consider two types of equilibria, the classical Stackelberg and Nash.

\subsubsection{Stackelberg Equilibrium}
In the Stackelberg equilibrium, there exists a leader, and a follower and the leader can enforce his strategy on the follower. It is assumed that the follower plays rationally, in the sense that he plays the best response with respect to the strategy of the leader. Thus, the leader chooses the strategy which maximizes his payoff given the best response of the follower. 

\begin{definition}[\cite{Basar_1995} \hspace{-0.1cm}] \label{def:Stackelberg}
A strategy pair $(i^*, j^*) \in \Gamma^1 \times \Gamma^2$ is a Stackelberg equilibrium of a bimatrix game $A$, $B$, with P1 as leader, if:
\begin{align*}
i^*  = \arg\max_{i \in \Gamma^1} \min_{j \in R(i)} a_{i,j}  \nonumber \,,  \quad j^* = \arg\max_{j\in \Gamma^2} b_{i^*,j}\,, 
\end{align*}
where $R(i) := \arg\max_{j\in \Gamma^2} b_{i,j}$ is the best response of P2 to the strategy $i \in \Gamma^1$ of P1.
\end{definition}

Note that the best trajectories of both the leader and the follower may be non-unique. For the leader all optimal trajectories lead to the same pay-off; the leader can, therefore, pick anyone among them and announce it to the follower \cite{Basar_1995}. The follower can then select any one of the (possibly many) best responses to this announced trajectory of the leader.

\subsubsection{Nash Equilibrium}
The Nash equilibrium is a trajectory pair, such that there is no incentive for either of the players to deviate unless the other player does. 
\begin{definition}[\cite{Basar_1995} \hspace{-0.1cm}] \label{def:nash}
A strategy pair $(i^*, j^*) \in \Gamma^1 \times \Gamma^2$ is a Nash equilibrium of a bimatrix game $A$, $B$ if the following inequalities are fulfilled:
\begin{align*}
a_{i^*,j^*} &\geq a_{i,j^*} \quad \forall i \in \Gamma^1\,, \\
b_{i^*,j^*} &\geq b_{i^*,j} \quad \forall j \in \Gamma^2 \,.
\end{align*}
\end{definition}
We note that, unlike zero sum games, existence of Nash pure strategy equilibria (as Definition 4) is not guaranteed for our bimatrix games, even though the payoff matrices are bounded; we address the issue of existence of pure strategy Nash equilibria below. Moreover, if Nash equilibria do exists, there may be multiple of them with potentially different payoffs for the two players \cite{Basar_1995}. To distinguish between them we use the notion of ``betterness" following terminology of \cite{Basar_1995}. 
\begin{definition}[\cite{Basar_1995} \hspace{-0.1cm}] \label{def:better}
A strategy pair $(i^1,j^1)$ is said to be better than another strategy pair $(i^2,j^2)$, if $a_{i^1,j^1} \geq a_{i^2,j^2}$, $b_{i^1,j^1} \geq b_{i^2,j^2}$, and at least one of the inequalities is strict.
\end{definition}

If neither equilibrium is better than the other we then say that the two equilibria are incomparable. The existence of better equilibria forms the basis for a non-cooperative equilibrium consensus \cite{Basar_1995}.

\subsubsection{Feasibility Assumptions}

Since feasibility of a trajectory pair is fundamental for the game, but the constraints are only considered in the payoff, we investigate assumptions under which a feasible trajectory pair is obtained. 

\edef\oldassumption{\the\numexpr\value{assumption}+1}
\setcounter{assumption}{0}
\renewcommand{\theassumption}{\oldassumption.\alph{assumption}}
\begin{assumption} 
There exists a feasible trajectory pair $(i,j)$ according to Definition \ref{def:feas}.\label{ass:Feas_a}
\end{assumption}

To establish a connection between the sequential and the cooperative game we need a more restrictive assumption involving the payoff of the leader (P1) in the sequential game \eqref{eq:SeqRacingGame}. 
\begin{assumption}
There exists a trajectory pair $(i,j)$, such that $\max_{i \in \Gamma^1} a_i > \kappa$ and $\max_{j\in \Gamma^2} \min_{i' \in \arg\max_{i \in \Gamma_1}} b_{i',j} > \lambda$.\label{ass:Feas_b} 
\end{assumption}
\let\theassumption\origtheassumption

Note that, the payoff of the follower is identical in both the cooperative game \eqref{eq:racingGame} and the sequential game \eqref{eq:SeqRacingGame}, thus the assumption is well defined for both games. One can easily verify that Assumption \ref{ass:Feas_b} is more restrictive, since any cooperative game which fulfills Assumption \ref{ass:Feas_b} also fulfills Assumption \ref{ass:Feas_a}, however, the opposite does not necessarily hold. We note that Assumption \ref{ass:Feas_b} is not always fulfilled in practice, but the assumption is necessary for our technical results.

\subsection{Equilibrium of the sequential Game}

When investigating the sequential game \eqref{eq:SeqRacingGame} in more detail one can see that it is not really a game since the decisions of the follower do not influence the decisions of the leader. Thus the optimal trajectory pair can simply be computed by a sequential maximization approach, where the leader selects an optimal trajectory without considering the follower. Similar to the Stackelberg equilibrium this trajectory is then announced to the follower, who selects a best response. This can be summarized in the following two-step sequential maximization, which determines the sequential optimal trajectory pair $(i^s, j^s)$,
\begin{align}
i^s =  \arg\max_{i\in \Gamma^1} a_i \,, \quad
j^s =  \arg\max_{j\in \Gamma^2} b_{i^s,j} \,.  \label{eq:seqMax}  
\end{align}

One can show that trajectory pairs $(i^s, j^s)$ that satisfy \eqref{eq:seqMax} are both Stackelberg and Nash equilibria of the sequential game \eqref{eq:SeqRacingGame}; the proof is similar to that of Theorem \ref{theo:Link} below and is omitted in the interest of space. Moreover, if Assumption \ref{ass:Feas_b} holds, these trajectory pairs are also feasible; this follows directly from Assumption \ref{ass:Feas_b}, which guarantees that the follower has a feasible response to the announced optimal trajectory of the leader. Once again multiple optimal trajectories for the leader do not pose a problem since their payoffs are the same and the chosen trajectory is announced to the follower.

In the racing situation in Fig. \ref{fig:raceSit} Assumption \ref{ass:Feas_b} is fulfilled and the sequential optimal trajectory pair is $(2,1)$, where the car in front goes straight and the car in the back avoids the car in front by going left. It is easy to verify that this trajectory pair is indeed a Nash and a Stackelberg equilibrium for \eqref{eq:ABraceSitSeq}.

\subsection{Equilibrium of the cooperative game}

Let us now focus on the more complex cooperative game \eqref{eq:racingGame} which does not exhibit the sequential structure of \eqref{eq:SeqRacingGame}. When investigating the racing situation in Fig. \ref{fig:raceSit} with the corresponding payoff matrices \eqref{eq:ABraceSit}, we see that the Stackelberg equilibrium of the game is the trajectory pair $(2,1)$, which is identical to the sequential optimal trajectory pair. Furthermore, we can also verify that the sequential optimal trajectory pair is also a Nash equilibrium. Note that this is not the only Nash equilibrium, however: the trajectory pair $(1,2)$ is also a Nash equilibrium. 

A natural question that arises is how can the two cars choose one of the two Nash equilibria. This is generally a difficult question in game theory. In our example, if each player assumes that the equilibrium that maximizes his own progress will be chosen, this would result in the trajectory pair (2,2) which is not a Nash equilibrium and results in a collision\footnote{Note that in this example the best-response dynamics as proposed in \cite{Williams2017} would cycle between trajectory pair $(2,2)$ and $(1,1)$.}. This is because none of the Nash equilibria is better than the other. To resolve this issue ahead of the game the players can agree on rules about how to pick a Nash equilibrium in the case of multiple equilibria; we call those ``rules of the road". Here we use a rule that says that if there are multiple Nash equilibria, the equilibrium with the largest payoff for the leader (P1) is chosen. In the racing situation in Fig. \ref{fig:raceSit}, the Nash equilibrium which fulfills the rules of the road is the trajectory pair $(2,1)$. 

Note that in this case all the different equilibrium concepts, as well as the sequential game, lead to the same \emph{feasible} trajectory pair. Indeed one can show that under Assumption \ref{ass:Feas_b} this property holds generally.

\begin{theorem}\label{theo:Link} $ $
\begin{itemize}
    \item[(a)] If Assumption \ref{ass:Feas_a} holds, all Stackelberg equilibria are feasible, there exists at least one feasible Nash equilibrium and all feasible Nash equilibria are better than all infeasible Nash equilibria.
    \item[(b)] Let $\Pi_\text{s}$ denote the set of all Nash/Stackelberg equilibria of the sequential game, $\Pi_\text{st}$ the set of all Stackelberg equilibria of the cooperative game, $\Pi_\text{n,RoR}$ the set of all Nash equilibria of the cooperative game that fulfill the rules of the road, and $\Pi_\text{n}$ the set of all Nash equilibria of the cooperative game. If Assumption \ref{ass:Feas_b} holds, then
\begin{align*}
    \emptyset \neq \Pi_\text{s} = \Pi_\text{st} = \Pi_\text{n,RoR} \subseteq \Pi_\text{n}\,.
\end{align*}
\end{itemize}
\end{theorem}

\begin{IEEEproof}
We start by noting that due to the payoff structure in the cooperative game \eqref{eq:racingGame}, specifically the symmetry of the collision constraints, it holds that 
\begin{align} \label{eq:sym}
a_{i,j} = \lambda \Leftrightarrow b_{i,j} = \lambda \,.
\end{align}

Part (a): Under Assumption \ref{ass:Feas_a} there exists $i \in \Gamma^1$ such that $b_{i,j} > \lambda$ for all $j \in R(i)$, hence by \eqref{eq:sym} $\max_{i \in \Gamma^1} \min_{j \in R(i)} a_{i,j} > \lambda$ and therefore all Stackelberg equilibria are feasible. 

We now show that the Stackelberg equilibrium is a Nash equilibrium, which by virtue of the previous point shows that a feasible Nash equilibrium exists. Note that $\max_{i \in \Gamma^1} \min_{j \in R(i)} a_{i,j} = \max_{i \in \Gamma^1} \max_{j \in \Gamma^2} a_{i,j}$ due to the payoff structure. This observation together with P2 playing best response shows that the Nash inequalities are fulfilled. 

Finally, note that for any infeasible Nash equilibrium at least one of the two players is not using their feasible trajectory, as otherwise, the Nash would be feasible. Moreover, none of the two players violates the track constraint; if they did, they could improve their payoff from $\kappa$ to (at least) $\lambda$ by switching to another trajectory, contradicting Nash. Hence neither is using their trajectory from the feasible pair, as if one was the other could unilaterally improve their payoff by also switching to their feasible trajectory, contradicting Nash again. Therefore an infeasible Nash equilibrium only exists if both players have a payoff of $\lambda$ but for both players, any unilateral change would result in a payoff of $\lambda$ or $\kappa$. However, any feasible Nash equilibrium is better than such an infeasible Nash equilibrium, as in the feasible case both players receive the progress payoff, which is by definition greater than $\lambda$.

Part (b): If Assumption \ref{ass:Feas_b} holds then $b_{i^s,j} > \lambda$ for all $j \in R(i^s)$, which again implies by \eqref{eq:sym} that $\min_{j\in R(i^s)} a_{i^s,j} = p(x^1_{i^s}(N))$. This also shows that the sequential optimal trajectory pair is feasible and the set of trajectory pairs is non-empty. 

Furthermore, the sequential optimal trajectory pair $(i^s,j^s)$ is identical to the Stackelberg equilibrium since by the above observation $\max_{i \in \Gamma^1} \min_{j \in R(i)} a_{i,j} = \max_{i \in \Gamma^1}  a_i$. Thus, we can reformulate the Stackelberg equilibrium as,
\begin{align*}
i^* &= \arg \max_{i \in \Gamma^1} \min_{j \in R(i)} a_{i,j} = \arg \max_{i \in \Gamma^1}  a_i\,,\\
j^* &= \arg\max_{j\in \Gamma^2} b_{i^*,j}
\end{align*}
which is identical to the sequential maximization approach \eqref{eq:seqMax} and proves the second equality. 

Finally, we show that the sequential optimal trajectory pair $(i^s,j^s)$ is a Nash equilibrium that fulfills the rules of the road. We know that $a_{i^s,j^s}  = p(x^1_{i^s}(N))$, leading to
\begin{align*}
&a_{i^s,j^s}  = p(x^1_{i^s} (N))  \geq a_{i,j^s} \quad \forall i\in \Gamma^1\,,\\
&b_{i^s,j^s} = b_{i^s,R(i^s)}  \geq b_{i^s,j}  \quad \forall j\in \Gamma^2\,.
\end{align*}
The first inequality holds since $p(x^1_{i^s}(N))$ is the best possible payoff for P1, the second because P2 plays best response. Therefore, the trajectory pair $(i^s,j^s)$ is a Nash equilibrium. Furthermore, the trajectory pair $(i^s,j^s)$ is the Nash equilibrium which fulfills the rules of the road, as it yields to the best possible payoff for P1, which shows the last two relationships. 
\end{IEEEproof}

\begin{figure}[ht]
\centering
\includegraphics[width = 0.475\textwidth]{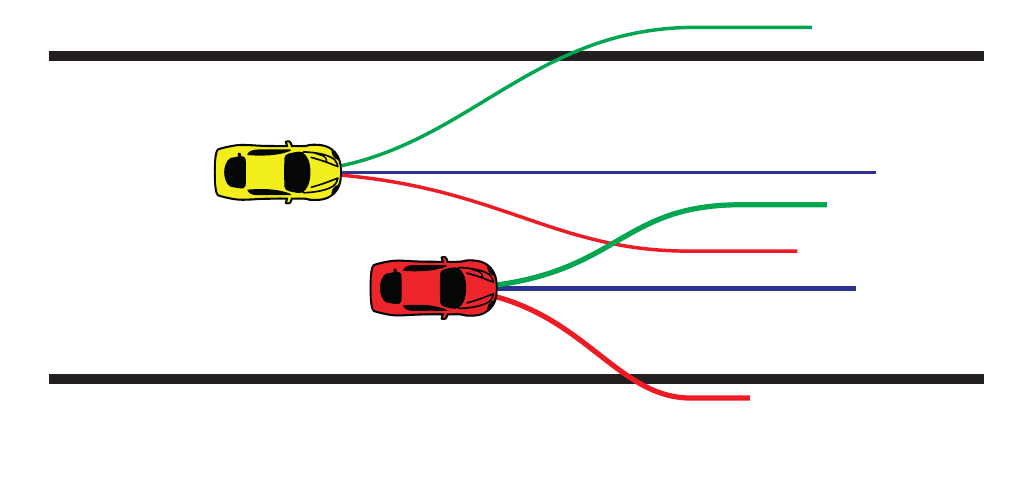}
    \caption{A racing game with an infeasible Nash equilibrium. Each player has three trajectories green, blue, red from top to bottom, and the trajectories are numbered from top (1) to bottom (3). The green trajectory (1) of P1 (red car) collides both with the blue (2) and the red (3) trajectory of P2 (yellow car) and the blue trajectory (2) of P1 also collides with the red trajectory (3) of P2. }\label{fig:InfeasNash}
\end{figure}

To illustrate the proof of Theorem \ref{theo:Link} (a), consider the situation depicted in Fig. \ref{fig:InfeasNash} with the payoff matrices in \eqref{eq:InfeasNash}. In this case, two Nash equilibria exist, $(2,2)$ which is feasible with a payoff of $(0.87,0.89)$, and $(1,3)$ which is infeasible with a payoff of $(-1,-1)$. The feasible Nash is better than the infeasible one and is thus preferable for both players.

\begin{align} \label{eq:InfeasNash}
A = \begin{bmatrix}
0.84 & -1 & -1\\
0.87 & 0.87 & -1 \\
-10 & -10 & -10
\end{bmatrix} \, B = \begin{bmatrix}
-10 & -1 &-1\\
-10 & 0.89 & -1 \\
-10 & 0.81 & 0.81
\end{bmatrix}.
\end{align}

Thanks to Theorem \ref{theo:Link}, if Assumption \ref{ass:Feas_b} holds we can use the sequential maximization approach \eqref{eq:seqMax} to compute Nash and Stackelberg equilibria. This can offer significant computational advantages since it is not necessary to generate all entries of $A$ and $B$. Moreover, collision checks are only necessary between the optimal trajectory of the leader and the possible trajectories of the follower; this dramatically reduces the computational burden as we will discuss in Section \ref{sec:SimResults}. We note that the only structure necessary for Theorem \ref{theo:Link} is the fact that the decisions of the players are only linked through the constraints as the progress payoff solely depends on the players own action. Thus, Theorem \ref{theo:Link} would still hold if one uses objective functions other than progress, e.g., fuel efficiency, or even in other applications that enjoy a similar structure. This suggests that the theorem may have wider applicability, as many robotic multi-agent problems with collision constraints can be formulated such that they have the same structure.

Note that if Assumption \ref{ass:Feas_b} is not fulfilled, in the cooperative game the leader will choose a trajectory which allows the follower to make a choice such that the resulting trajectory pair is feasible (as long as Assumption \ref{ass:Feas_a} holds). On the other hand, in the sequential game the leader does not consider the follower and a trajectory pair that satisfies \eqref{eq:seqMax} may be not feasible even if Assumption \ref{ass:Feas_a} is satisfied.

\subsection{Equilibrium of the blocking game}

Finally let us investigate the blocking game. To formalize the discussion, we start by defining a blocking trajectory pair. Let $(i^{cg}, j^{cg})$ denote a Stackelberg equilibrium of the cooperative game \eqref{eq:racingGame}. Note that if there are multiple Stackelberg equilibria, all these trajectory pairs have the same progress payoff. Therefore, all the progress based arguments still hold in the case of multiple equilibria.

\begin{definition} \label{def:Blocking}
A trajectory pair $(i,j)$ is blocking if the following properties hold: $(i)$ the Stackelberg equilibrium of the cooperative game $(i^{cg}, j^{cg})$ is such that P1 gets overtaken at the end of the horizon $p(x_{i^{cg}}^{1}(N)) < p(x_{j^{cg}}^2(N))$, $(ii)$ the pair $(i,j)$ is such that P1 does not get overtaken $p(x_{i}^{1}(N)) > p(x_{j}^{2}(N))$, $(iii)$ the pair $(i,j)$ is feasible $a_{i,j} > \lambda$ and $b_{i,j} > \lambda$, $(iv)$ for all $j^c \in \Gamma^2$ such that $p(x_{j^c}^{2}(N))  > p(x_{j}^{2}(N))$, it holds that $b_{i,j^c} \leq \lambda$. 
\end{definition}

A blocking trajectory pair corresponds to a collision-free trajectory pair where P1 is ahead at the end of the horizon. P1 achieves this by choosing a trajectory, such that any trajectory of P2 that would achieve a larger progress collides with this trajectory of P1. The last condition also implies that P2 plays a best response, as any other trajectory would lead to a smaller payoff. Among all blocking trajectory pairs fulfilling properties $(i)$-$(iv)$, let $(i^b,j^b)$ denote the one which achieves the largest progress for the leader.

When investigating the racing situation in Fig. \ref{fig:blocking} with the payoff matrices \eqref{eq:BlockingEx}, we can see that the additional payoff term promotes blocking behavior since the Stackelberg equilibrium $(2,1)$ of the game is a blocking trajectory where P1 drives a curve to prevent P2 from overtaking at the end of the horizon. In contrast, the (unique) Nash equilibrium is $(3,2)$, where P1 drives straight to maximize progress, but P2 can overtake P1 at the end of the horizon. Interestingly, this is the same trajectory pair which would be optimal if the cooperative game payoff would be used. Similar to the cooperative game, the observations from the simple racing situation can again be generalized.

\begin{theorem} \label{theo:StackelbergBlockingGame}
Under Assumption \ref{ass:Feas_a}:
\begin{itemize}
\item[(a)] If $p(x_{i^{cg}}^{1}(N)) > p(x_{j^{cg}}^2(N))$, or if $p(x_{i^{cg}}^{1}(N)) \leq p(x_{j^{cg}}^2(N))$ but there does not exists a blocking trajectory pair, then the Stackelberg equilibrium $(i^{cg}, j^{cg})$ of the cooperative game \eqref{eq:racingGame} is identical to the Stackelberg equilibrium of the blocking game \eqref{eq:blockingGame}. If there exists a blocking trajectory pair and if $p(x_{i^{cg}}^{1}(N)) \leq p(x_{i^b}^1(N)) + w$, then $(i^b,j^b)$ is the Stackelberg equilibrium of the blocking game \eqref{eq:blockingGame}. 

\item[(b)] The Stackelberg equilibrium $(i^{cg}, j^{cg})$ of the cooperative game \eqref{eq:racingGame} is a Nash equilibrium of the blocking game \eqref{eq:blockingGame}. And if a blocking trajectory $(i^b,j^b)$ is a Nash equilibrium of the cooperative game \eqref{eq:racingGame} it is a Nash equilibrium of the blocking game \eqref{eq:blockingGame}.

\end{itemize}

\end{theorem}

\begin{IEEEproof}
We start by noting that by Assumption \ref{ass:Feas_a} and Theorem \ref{theo:Link} we know that $(i^{cg}, j^{cg})$ is feasible, and since the constraints are the same in the cooperative game and the blocking game, we know that this trajectory pair is also feasible for the blocking game.

Part (a): Given this observation we know that if $p(x_{i^{cg}}^{1}(N)) > p(x_{j^{cg}}^2(N))$, it follows that the payoffs for the blocking game of this trajectory pair are, $a_{i^{cg}, j^{cg}} = p(x_{i^{cg}}^{1}(N)) + w$ and $b_{i^{cg}, j^{cg}} = p(x_{j^{cg}}^{2}(N))$. We now show that $j^{cg}$ is the best response of P2 to $i^{cg}$ given the payoff of the blocking game \eqref{eq:blockingGame}. We know that $j^{cg}$ maximizes the progress given $i^{cg}$ due to the formulation of the cooperative game \eqref{eq:racingGame}, and since the additional payoff only depends on the progress, we know that no other trajectory would get the additional reward, thus $j^{cg}$ is the best response. Second, $i^{cg}$ maximizes the payoff for P1, therefore, $(i^{cg}, j^{cg})$ is also the Stackelberg equilibrium of the blocking game.

Second, if $p(x_{i^{cg}}^{1}(N)) \leq p(x_{j^{cg}}^2(N))$ and there exists no blocking trajectory pair $(i^b,j^b)$, or in other words P1 cannot avoid the overtaking maneuver, the payoff of P1 does not depend on the blocking payoff, and is therefore maximized by $i^{cg}$. As $j^{cg}$ leads to the largest possible progress for P2 and by definition gets the additional reward $w$, $j^{cg}$ is also the best response for the blocking game \eqref{eq:blockingGame}.

Third, let us assume there exists a blocking trajectory pair $(i^b,j^b)$. The pair will be chosen if it leads to the largest payoff for the two players. For P1 this is the case if $p(x_{i^b}^{1}(N)) + w \geq a_{i,j}$ for all $i$ and $j$, since $i^b$ is the best possible blocking trajectory other blocking strategies have a lower payoff. Thus, the inequality holds if the payoff $a_{i^b,j^b}$ is larger than the largest payoff of a feasible non-blocking trajectory pair, which is $p(x_{i^{cg}}^{1}(N))$. Thus P1 chooses $i^b$ if $p(x_{i^{cg}}^{1}(N)) \leq p(x_{i^b}^1(N)) + w$ and by definition $j^b$ is the optimal response for P2. Thus $(i^b,j^b)$ is a Stackelberg equilibrium.

Part (b): To show that $(i^{cg}, j^{cg})$ is a Nash equilibrium of the blocking game \eqref{eq:blockingGame}, two cases are possible. First, $p(x_{i^{cg}}^{1}(N)) > p(x_{j^{cg}}^2(N))$, where $a_{i^{cg},j^{cg}} = p(x_{i^{cg}}^{1}(N)) + w$ and $b_{i^{cg},j^{cg}} = p(x_{j^{cg}}^2(N))$, and second,  $p(x_{i^{cg}}^{1}(N)) \leq p(x_{j^{cg}}^2(N))$, where $a_{i^{cg},j^{cg}} = p(x_{i^{cg}}^{1}(N))$ and $b_{i^{cg},j^{cg}} = p(x_{j^{cg}}^2(N)) + w$. In the first case, because $j^{cg}$ is the best response in terms of progress we know that if $j^{cg}$ does not get the additional reward, no other trajectory gets the reward. Therefore, it holds that for all $j\in \Gamma^2$, $ b_{i^{cg},j^{cg}} \geq b_{i^{cg},j}$. As a consequence for P2 the trajectory pair $(i^{cg},j^{cg})$ fulfills the Nash inequality. For P1 the Nash inequality $a_{i^{cg},j^{cg}} \geq a_{i,j^{cg}}$ is fulfilled since it is the largest feasible payoff for P1. In the second case, a similar argument holds, which is omitted in the interest of space. Thus, $(i^{cg},j^{cg})$ is a Nash equilibrium of the blocking game \eqref{eq:blockingGame} even though the additional payoff is not considered in the computation of the trajectory pair.

For the second part, we now assume that the Nash equilibrium is a blocking trajectory pair $(i^b,j^b)$. The Nash inequality for P2 $b_{i^b,j^b} \geq b_{i^b,j}$ is always fulfilled by Definition \ref{def:Blocking}, since no action of P2 can get the additional payoff $w$ and $j^b$ is the best response of P2. However, the Nash inequality for P1,
\begin{align*}
p(x_{i^b}^1(N)) + w \geq a_{i,j^b} \quad \forall i\in \Gamma^1\,,
\end{align*}
is only fulfilled, if for all $i\in \Gamma^1$ such that $p(x_{i}^{1}(N)) > p(x_{i^b}^{1}(N))$, it holds that $a_{i,j^b} \leq \lambda$. This condition, however, is identical to the one required by the Nash equilibrium of the cooperative game \eqref{eq:racingGame}. Thus, for $(i^b,j^b)$ to be a Nash equilibrium of the blocking game \eqref{eq:blockingGame}, it also needs to be a Nash equilibrium of the cooperative game \eqref{eq:racingGame}.

\end{IEEEproof}

This illustrates, that if the goal is to encourage blocking the Stackelberg equilibrium is the appropriate equilibrium concept. Furthermore, this also shows that blocking needs an asymmetric information pattern, which is logical as the leader needs to be sure that he can enforce the blocking trajectory on the follower. 

One can recognize that the blocking game subsumes the cooperative game, more precisely for $w = 0$ the two games are identical. This poses the question for which values of $w$ the equilibrium changes. This can be answered for the Stackelberg equilibrium using part (a) of Theorem \ref{theo:StackelbergBlockingGame}. For a game instance where a blocking trajectory pair $(i^b,j^b)$ exists, there is a discrete change in the optimal trajectory pair when increasing $w$. This switch occurs the moment $p(x_{i^b}^1(N)) + w$ becomes larger that $p(x_{i^{cg}}^{1}(N))$; in \eqref{eq:BlockingEx} the switch occurs if $w \geq 0.03$. From a racing point of view, it is best to choose $w$ large enough such that whenever possible the blocking trajectory is selected.

\section{Moving Horizon Games} \label{sec:MovingHorizonGames}

The racing games presented in Section \ref{sec:GameForm} are finite horizon games. To introduce feedback we propose to play the games in a moving horizon fashion, solving the game with the full horizon but only applying the first input, then reformulating the game based on the current state measurement and solving again. Similar to \cite{Cruz2002} and \cite{Virtanen2006}, this allows to approximately solve an infinite horizon game, which would not be tractable using dynamic programming, by solving a series of finite horizon games. 

An inherent problem of moving horizon approaches is that they, in general, do not guarantee feasibility of closed-loop operation. This problem can be tackled by using terminal set constraints which ensure recursive feasibility by guaranteeing that the terminal state is within an invariant set \cite{Mayne2000}. Alternatively one can use soft constraints, rendering the problem feasible at all times \cite{Kerrigan2000}.

\subsection{Terminal viability constraints}
If the goal is to guarantee feasibility of the racing games at all times when played in a receding horizon fashion, one can use tools from viability theory. Viability theory is concerned with finding the set of initial conditions for which there exists a viable solution, i.e., a solution to a difference inclusion which forever remains within a given closed constraint set $K$ \cite{Aubin2009}. For a difference inclusion $x_{k+1} \in F(x_k)$ this set of initial conditions is called the \emph{viability kernel} (denoted by $\text{Viab}_F(K)$) and can be numerically approximated by the viability kernel algorithm, see \cite{saintPierre94} for more details.  

In \cite{Liniger_2017} viability theory was already used to guarantee recursive feasibility with respect to track constraints in an autonomous racing set up. Here we extend this approach to also consider viability with respect to collisions between the cars, relying on the simple observation that, once in the viability kernel there exists an input which keeps the system state in the viability kernel \cite{Cardaliaguet99}. We distinguish the sequential game from the cooperative and blocking games since the former requires the more restrictive Assumption \ref{ass:Feas_b} to hold whereas for the latter two Assumption \ref{ass:Feas_a} suffices to guarantee feasibility.

For the sequential game the aim is to fulfill Assumption \ref{ass:Feas_b} recursively this causes technical difficulties due to discontinuities that arise in the dynamics which need to be tackled appropriately. One discontinuity arises because the leader and follower may switch roles from one solution of the game to the next. This problem can be resolved using a hybrid state in combination with a spatial regularization to prevent unwanted switching behavior, see \cite{Johansson1999}. Note that the system is now hybrid requiring a hybrid viability kernel algorithm \cite{Margellos2013}. Second, the optimal viable policy of the leader leads to discontinuities which can be tackled by appropriately smoothing the policy of the leader. Given this continuous dynamical model incorporating the sequential maximization approach, the viability kernel with respect to the track and collision constraints can be computed. If the game starts with the car states in this viability kernel then the follower always has a viable trajectory if the leader plays his optimal viable trajectory.

The cooperative and blocking game, require a conceptually simpler viability kernel since guaranteeing Assumption \ref{ass:Feas_a} recursively can be assured if there always exists a \emph{cooperative} collision avoiding trajectory pair. Therefore, both cars need a viable trajectory with respect to track and collision constraints something that can be ensured by computing a standard viability computation in the product space.

For all three games, we can theoretically compute an appropriate viability kernel, which can be imposed as a terminal constraint to guarantee recursive feasibility. This can be achieved by penalizing trajectory pairs for which the terminal state is not in the viability kernel by a payoff $\lambda$. If the viability kernel is empty safe racing is not possible, as the cars are doomed to either collide or leave the track eventually; the same is true if the viability kernel is not empty but does not contain the initial condition of the game.

Unfortunately, practically it is not possible to compute the described viability kernels as they require one to consider the state of both players simultaneously. As this joint state is ten-dimensional the computation is not tractable with the standard, grid based viability kernel algorithm. From a practical point of view, it is possible to use the viability kernel with respect to the track constraints, as proposed in \cite{Liniger_2017}, as a first approximation. This computation is tractable \cite{Liniger_2017} and it is easy to see that the product of the track viability kernels for the two cars is a superset of all racing viability kernels discussed above; ensuring that the system remains in this product set therefore provides a necessary condition for the recursive feasibility of the games. We adopt this approach in our simulation study below.

\subsection{Soft Constraints}
Alternatively one can use soft constraints that allow formulating a game which always finds an optimal trajectory. The idea is to relax the constraints and consider them by reducing the payoff depending on how much the constraints are violated. The advantage compared to the above discussed terminal constraints is that no viability kernels need to be computed. However, the resulting trajectory is not always feasible with respect to the constraints.

The idea is similar to the soft constraints in MPC \cite{Kerrigan2000}, however, the implementation in the game set-up is quite different. In our racing games, we use soft constraints for the collision constraints only. In a first step, given a trajectory pair we can compute the smallest relaxation of the constraint necessary, which we call slack multipliers of a trajectory pair $S_{i,j} = \sum_{k=1}^N \max(-d(x^1_i(k),x^2_j(k)),0)$. Given our collision constraints, this corresponds to how much the two cars penetrate each other and therefore $S_{i,j}$ is always positive. Second, instead of penalizing a colliding trajectory pair with a payoff of $\lambda$, the soft constrained game computes the payoff neglecting the collision constraints and then reduces this payoff by $\sigma S_{i,j}$, where $\sigma$ is a positive weight. Thus, if $\sigma$ is chosen large enough, the players in the game have still no benefit from violating the collision constraints, and therefore the optimal trajectory pair does not change. This is often called an exact soft constraint reformulation, and it can be shown that there always exist such a $\sigma$ for all racing games. The main reason this holds is that there is only a finite number of trajectories and all payoffs as well as the slack multipliers $S_{i,j}$ are bounded. On the other hand, if there is no feasible trajectory the soft constrained game does still find a reasonable trajectory pair by trading off violation of the collision constraint and payoff depending on the value of $\sigma$.

\section{Simulation Results} \label{sec:SimResults}

\subsection{Simulation set-up}
Our simulation study replicates the miniature race car set-up, hosted at the automatic control laboratory at ETH Z\"urich \cite{Liniger_2014}. The set-up consists of miniature race cars measuring 12$\times$5\,cm, driving around an approximately 18\,m long track (\url{https://youtu.be/RlZdMojOni0}). To simulate the cars we use a dynamical bicycle model with nonlinear tire forces as described in \cite{Liniger_2014}. 

\begin{figure}[ht]
\centering
\includegraphics[width = 0.45\textwidth]{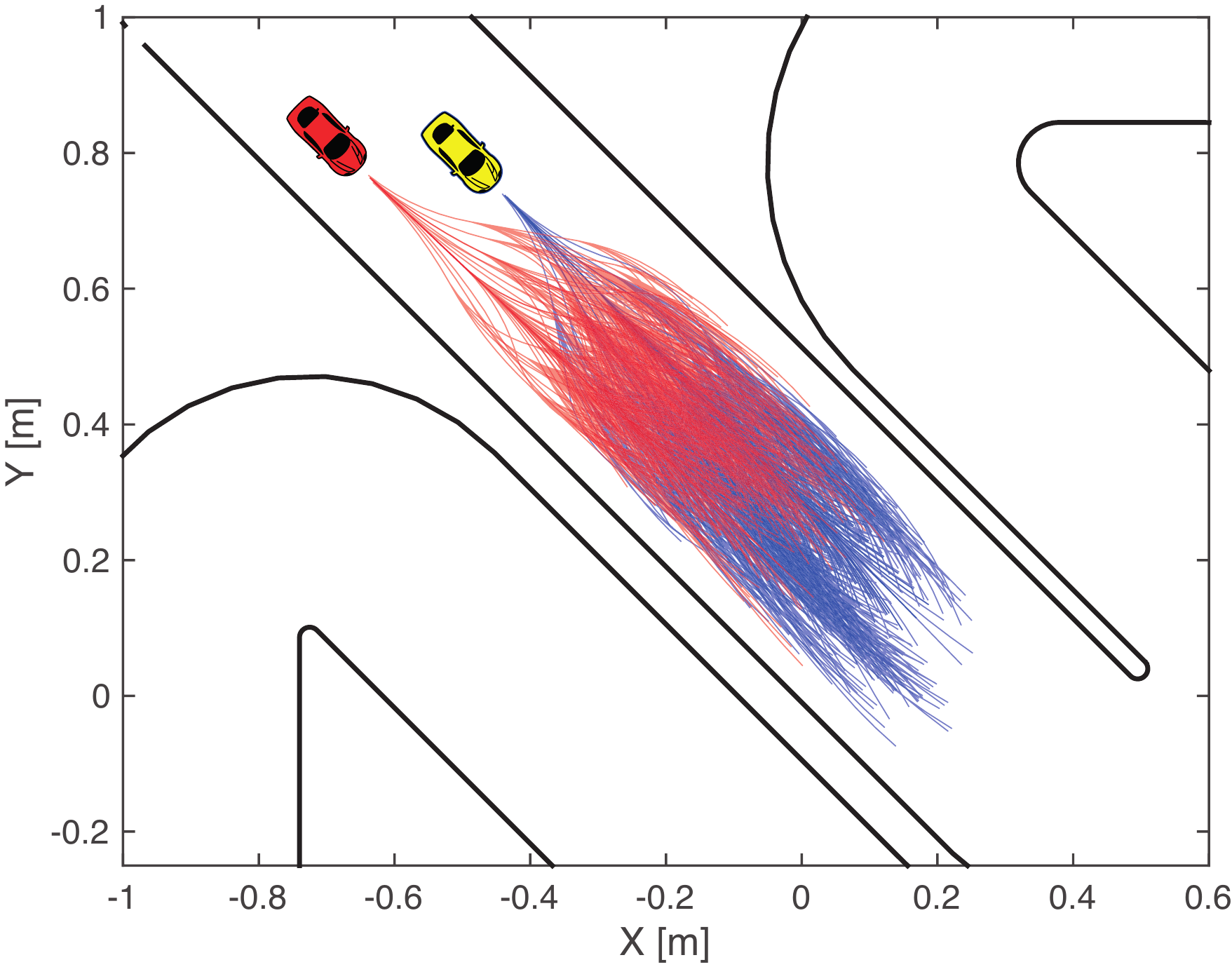}
    \caption{Possilbe trajectories of the path planning model for two cars using the viability-based pruning presented in \cite{Liniger_2017}}\label{fig:raceSitReal}
\end{figure}

To achieve high-performance driving we use a two-level hierarchical controller following the structure of \cite{Liniger_2017}. The upper level of the controller generates a trajectory pair based on one of the three racing games and the path planning model studied in the experimental set-up of \cite{Liniger_2017}. The path planning model makes use of the so-called constant velocity segments to generate a tree of alternative paths for each car. To provide a longer prediction horizon at a manageable computation burden longer sampling times and fewer prediction steps are used. Following \cite{Liniger_2017}, for the path planning model we use number prediction steps $N = 3$, discretization time $T_{pp} = 0.16$\,s and number of constant velocities $N_m = 129$.  Fig. \ref{fig:raceSitReal} shows an example of trajectories resulting from the path planning model. At the lower level, we employ an MPC controller to track the trajectories corresponding to the optimal trajectory pair. The lower level MPC uses the linearized bicycle model, and also enforces track constraints; more details can be found in \cite{Liniger_2017}. Note that the two levels are coupled since the lower level needs a trajectory to track, which is computed by one of the racing games, and the higher level racing game cannot generate inputs to the car which allow driving at the limit.

Note that brute force implementation of the racing bimatrix game would, in this case, require forming two $129^3\times 129^3$ payoff matrices, clearly a formidable computational task. To reduce the computation we rely on viability-based pruning of the trajectories \cite{Liniger_2017}, which excludes trajectories that are doomed to leave the track either immediately or at a later stage. Thus, these pruning strategies are non-heuristic since the trajectories would have a payoff of $\lambda$ in the long run. Following \cite[Section 4]{Liniger_2017} we use either the viability or the discriminating kernel to prune the trajectories. The two kernels lead to comparable lap times but differences in driving style that, in the context of the present paper, result in interesting racing behaviors. More precisely, pruning based on the discriminating kernel results in a more conservative driving style, breaking earlier into curves but achieving higher exit velocities out of curves. Pruning based on the viability kernel, on the other hand, results in a more aggressive driving style. In the following, we will use \emph{PV} to denote the player using the viability kernel and \emph{PD} the player using the discriminating kernel to prune the trajectories.

To generate the bimatrix game, we assume that each car shares their possible trajectories with the opponent. Alternatively, it would also be possible to share the path planning model and the pruning strategy at the beginning of the race, and online only share the current state. For real racing, this can be extended by estimating the opponents constant velocity points as well as estimating the state based on sensors. The estimation of the constant velocity points is simplified by the fact that similar cars have similar constant velocities. Thus, each car only needs to adapt its own model to get an estimate of the opponent's model. Alternatively, for the cooperative game a modified best response method can be used which is guaranteed to converge to a Nash equilibrium \cite{Liniger_2018}, where only a sequence of best response trajectories needs to be shared to solve each game.

To further reduce the computation time, we compute collisions using a two-step collision detection: first, the distance between the center of the cars is checked, only if the distance is small enough for a collision to be possible in a second step a separating axis collision detection algorithm is used \cite{Gottschalk1996}. Finally, we compute progress of all remaining trajectories. The projection \eqref{eq:proj} is computed by first generating offline a piecewise affine approximation of the center line, comprising 488 pieces; online the projection can then be computed by finding the closest affine piece and taking an inner product. All the tasks except the path planning model run at a sampling time of 20\,ms allowing for accurate collision checks and fast feedback. To be able to react fast to the changes of the opponent player and to deal with model mismatch the game is also repeated every 20\,ms.

The last implementation detail is concerned with the soft constraints, any constraint violation of more than 1\,cm is treated as a collision between the cars. Note that thereby the problem is not always feasible and that in a case where no feasible solution can be found the car behind initiates an emergency braking maneuver, until the problem is feasible again. 

\subsection{Results}

The racing games were evaluated in a simulation study, where the cars are initiated at 500 different random initial positions along the center line of the track. To guarantee feasible initial conditions the cars are initiated at a low forward velocity of 0.5\,m/s, tangentially to the center line. To promote interesting racing situations we choose the initial position of the two cars in close proximity (0-20\,cm apart) and the physical parameters of both cars are identical, including engine power, braking power, and tire models. For each initial position, the simulation is run for 40\,s, which corresponds to approximately four laps. 

We compare three cases of the presented racing games. In the first case the optimal strategy pair is computed using the Sackelberg equilibrium of the sequential game \eqref{eq:SeqRacingGame}, in the second case that of the cooperative game and in the third case that of the blocking game \eqref{eq:blockingGame} with $w=100$. Recall that the solutions of the first two approaches are identical if Assumption \ref{ass:Feas_b} is fulfilled, but may differ otherwise. We also assume that the players play truthfully. This is reasonable since deviating from this trajectory leads to a lower payoff, either due to a collision or to worse progress. We note, however, that in the blocking game the leader could benefit from being untruthful. Simulation studies also showed that other methods such as random play do lead to a worse payoff, in some cases even for both players.

\begin{table}[ht]
\caption{Simulation study of the racing games}
\label{tab:sim}
\centering
\ra{1.3}
\begin{tabular}{@{}l c c c c @{}}\toprule
 & seq. game & coop. game & blocking game\\
 & (\emph{PV/PD}) &  (\emph{PV/PD}) &  (\emph{PV/PD})\\
 \midrule
\# of overtaking & 113       & 857           & 414 \\
maneuvers         & (4/109) & (502/355) & (231/183)    \\ 
\parbox[b][{0.6cm}][b]{2.3cm}{\# of runs with overtaking maneuvers}  & 106  & 309 & 220 \vspace{0.1cm} \\ 
collision probability  & $5.42\cdot 10^{-3}$  & $4.59\cdot 10^{-3}$ & $5.28\cdot 10^{-3}$  \\ 
\bottomrule
\end{tabular}
\end{table}  

Let us start by noting that in all cases we observe a substantial number of overtaking maneuvers (see Table \ref{tab:sim}), even though both cars have the same power, and only differ in the pruning strategies they employ. However, it is visible that the sequential game has the fewest overtaking maneuvers. The cooperative game has by far the most overtaking maneuvers, whereas the blocking game has approximately half as many overtaking maneuvers. It is also interesting to see that the cooperative and blocking game lead to a more even distribution of overtaking maneuvers between the different pruning strategies. This stands in contrast to the sequential game, where nearly all overtaking maneuvers are performed by the player using the discriminating kernel as a their pruning strategy (\emph{PD}). When looking at the probability of a collision, we can see that all three implementations have a low empirical collision probability of around $5\cdot 10^{-3}$, especially given the close initial proximity of the cars. The cooperative game has the lowest collision probability and the sequential game the highest. Note that some of the collisions can also be caused by the fact that the low-level MPC does not consider the collision constraints. Even though experimental and simulation results suggest that the low-level MPC can follow the trajectory very precisely, the games are sometimes so tight that even sub-millimeter tracking errors can result in a collision.

The main reason for the observed trends is first and foremost the fact that in the sequential game the leader does not consider the follower. This explains the fewer overtaking maneuvers, the higher collision probability, and the observation that \emph{PV} has nearly no overtaking maneuvers. Almost all overtaking maneuver need a certain cooperation between the cars, and this holds especially when \emph{PV} overtakes \emph{PD}. In the sequential game the leader does not consider the follower and ignores the risk of a collision. Due to the layout of the track, this is an advantage for \emph{PD}, which can prevent nearly all overtaking maneuvers, whereas \emph{PV} is not able to do so. In the cooperative and blocking games the leader adapts his strategy to accommodate the follower. This helps to prevent collisions, but at the same time allows for overtakes by the opposing car. In the blocking game the leader does consider the follower, and helps to prevent a collision, but at the same time actively tries to avoid overtaking maneuvers. This, in the end, results in a middle ground both in terms of overtaking and collision probability.
 
\begin{table}[ht]
\caption{Payoff of the racing games}
\label{tab:payoff}
\centering
\ra{1.3}
\begin{tabular}{@{}l c c c c @{}}\toprule
 & seq. game & coop. game & blocking game\\
 & (\emph{PV/PD}) &  (\emph{PV/PD}) &  (\emph{PV/PD})\\
 \midrule
 mean progress [m]  & 82.89  & 82.49 & 82.53 \vspace{0.1cm} \\ 
\# of stay  & 384       & 292           & 329 \\
ahead runs & (148/236) & (230/62) & (196/133)    \vspace{0.1cm}  \\ 
\# of wins  & (156/344)  & (412/88) & (307/193)  \\ 
\bottomrule
\end{tabular}
\end{table}  

Similar trends can also be seen when analyzing the overall payoff of the game. Table \ref{tab:payoff}, shows the mean progress of the different games, in how many runs the car starting ahead was able to also be ahead at the end and how many times which player won. The first two metrics are directly the ones we try to maximize in our games, and the third is, in a sense, the long term objective of the players. We can see that the sequential game achieves the largest progress of the three games, whereas the cooperative and the blocking game achieve very similar values. The sequential game achieves the largest progress because the leader does not adapt his strategy thus only optimizing the progress. In the other games, there are also more overtaking maneuvers which slow down the cars. The number of stay ahead runs, and winning runs show a similar trend to overtaking in Table \ref{tab:sim}. The sequential game is the best to stay ahead and win but does this by causing more and also deeper collisions. From the stay ahead and winning runs, it is even better to see how the additional term in the blocking game helps to be more competitive compared to the cooperative game.
 
To further highlight the difference between the three different game approaches, one of the 500 simulation runs is shown in a video which can be found at: \url{https://youtu.be/3Skl5qeFum8} for the sequential game, \url{https://youtu.be/Cp9OchB2S_M} for the cooperative game, and \url{https://youtu.be/Xxa8W9D3Z_A} for the blocking game. The videos emphasize the points discussed above, where the sequential game has zero overtaking maneuvers but there are collisions, in the blocking game there is only one overtaking maneuver, but in the cooperative game there are three overtaking maneuvers.

\begin{figure}[ht]
\centering
\includegraphics[width = 0.45\textwidth]{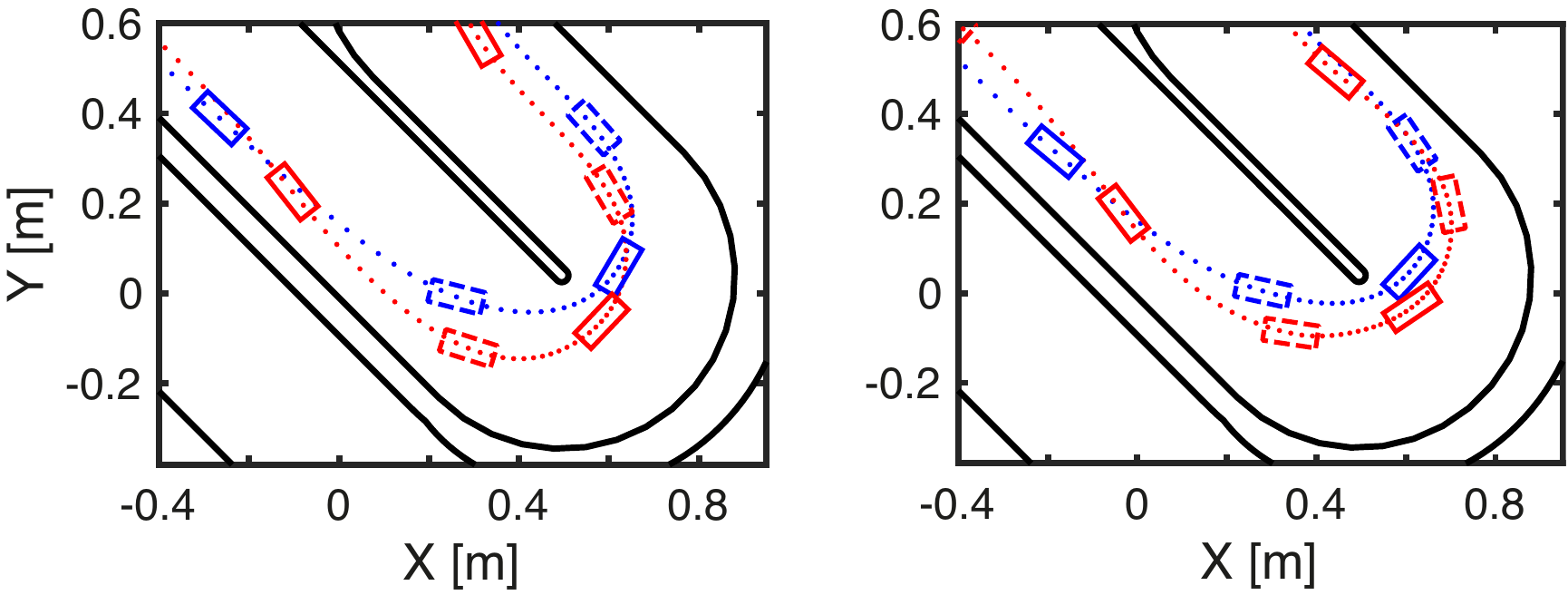}
    \caption{Time elapse of two overtaking maneuvers. Left side cooperative game, right blocking game. The blue car is \emph{PV} and the red car \emph{PD}, the cars are shown every 0.3\,s and every second time the car is dashed to highlight which cars are at the same time step. The driven trajectory also is also shown as a dotted line, where each dot is 0.02\,s apart. In both cases, the blue car can overtake the red car, but in the blocking game, the red car drives a more aggressive trajectory which tries to prevent the overtaking maneuver. }\label{fig:raceTimeLaps}
\end{figure}
 
\subsubsection{Remarks on difference in driving style and car dynamics}
The above simulation study investigated the case where one player \emph{PV} uses the viability and the other player \emph{PD} the discriminating kernel to prune their trajectories. Thus, the driving style of the two players is quite different, with \emph{PV} being ``aggressive" and \emph{PD} ``cautious." However, both players have the same model and identical cars. To further investigate the influence of different driving styles and differences in the cars, we also investigated the effect of a reduced motor power of one of the car as well as using the same approach to prune the trajectories for both players; the latter results in a race between identical cars and drivers. The simulation results show as expected that reducing the motor power of one car leads to more overtaking maneuvers. However, the influence is smaller than expected as the power constraint mainly influenced the maximal velocity, which is not too important in the used track layout. Using the same pruning strategy for both players drastically reduced the amount of overtaking maneuvers. However, even though the drivers and the cars are identical overtaking maneuvers still occur. In conclusion, we noticed that the change in the driving style induced by different pruning strategies seems to be more important than power differences between the cars.
\subsection{Computation time}
Simulations were carried out on a computer running Debian equipped with 16 GB of RAM and a 3.6 GHz Intel Xeon quad-core processor. When investigating the computation times reported in Table \ref{tab:compTime}, one can see that the computation times for the trajectory generation step is very similar for all three games, which is expected as the different cases only have a minor influence on the trajectory planning phase. However, the collision check for the sequential game is significantly faster than the collision checks for the other two games. 
\begin{table}[ht]
\caption{Computation time of the racing games}
\label{tab:compTime}
\centering
\ra{1.3}
\begin{tabular}{@{}l c c c c @{}}\toprule
 & seq. game & coop. game & blocking game\\
 \midrule
mean traj. generation [s] & 0.0027  & 0.0027 & 0.0027  \\ 
max traj. generation [s] & 0.0358  & 0.0314 & 0.0317  \\ 
mean collision check [s]  & 0.0012  & 0.3365 & 0.2545  \\ 
max collision check [s] & 0.0738 & 109.81 & 112.31 \\
\bottomrule
\end{tabular}
\end{table}  

The reason is that the cooperative and blocking game require checking collisions of all trajectories of both players with each other. Roughly speaking this computation grows quadratically in the number of trajectories available to the players and can result in up to $10^9$ collision checks. In the sequential game, when using the sequential maximization approach, collisions are only checked between the best trajectory of the leader and all trajectories of the follower. In this case, the number of collision checks grows only linearly in the number of trajectories available to the players. We can say that the sequential game is close to real-time feasible, whereas the cooperative and blocking game are significantly slower and would need radical changes to implement in real-time.

\section{Conclusion} \label{sec:Conclusion}
We presented a non-cooperative game approach for one-to-one racing of two autonomous cars. Three approaches to model the interaction between the two players were considered and formulated as bimatrix games. In the first two approaches the interaction is limited to collision constraints, and each player optimizes his progress. The third approach augments the cost function by rewarding staying ahead at the end of the horizon to emphasize blocking behavior. For online implementation the game is played in a moving horizon fashion and two methods were proposed to deal with the loss of feasibility in closed-loop. A simulation study shows that the proposed games can be used for competitive autonomous racing of two cars. The main observation is that the sequential game where the follower is completely neglected seems to be the most efficient blocking technique in closed-loop, but also comes with the highest risk of a collision. 

In future work, we plan to implement the games on our experimental set-up. To make this possible the computation time has to be substantially decreased. We envision accomplishing this by, among others, more efficient heuristic trajectory pruning techniques and exploiting the parallelizability of the method. Finally, we also plan to incorporate the collision constraints in the low-level MPC problem using the method proposed in \cite{Zhang_2018}.

\begin{appendices}

\section{Path planning model} \label{app:ppModel}
The path planning model is based on a dynamic bicycle model with nonlinear tire forces. The states of the system are the positions $X$ and $Y$, the heading angle $\varphi$, the longitudinal and lateral velocity $v_x$ and $v_y$, as well as the yaw rate $\omega$. The control inputs are the steering angle $\delta$ and the pulse width modulation duty cycle $d$ of the drive train motor. The equations of motion are
	\begin{align}
	\dot{X} &= v_x \cos(\varphi) - v_y \sin(\varphi)\,, \nonumber \\
	\dot{Y} &= v_x \sin(\varphi) + v_y \cos(\varphi)\,, \nonumber\\
	\dot{\varphi} &= \omega\,,\label{eq:odeModel}\\
	\dot{v}_x &= \frac{1}{m}\Big(F_{r,x}(v_x,d) - F_{f,y}(v_x,v_y,\omega,\delta) \sin{\delta} + m v_y \omega \Big)\,, \nonumber\\
	\dot{v}_y &= \frac{1}{m} \Big( F_{r,y}(v_x,v_y,\omega) + F_{f,y}(v_x,v_y,\omega,\delta) \cos{\delta} - m v_x \omega \Big) \,,\nonumber\\
	\dot{\omega} &= \frac{1}{I_z}\Big(F_{f,y}(v_x,v_y,\omega,\delta) l_f \cos{\delta} - F_{r,y}(v_x,v_y,\omega) l_r \Big)\,,\nonumber
	\end{align}
where $m$ is the mass of the vehicle, $I_z$ is the moment of inertia, and $l_r$ and $l_f$ are the distances from the CoG to the rear and the front wheel, respectively. The most important part are the forces modeling the tires, where $F_{r,x}(v_x,d)$ is the drivetrain force, $F_{r,y}(v_x,v_y,\omega)$ and $F_{f,y}(v_x,v_y,\omega,\delta)$ are the lateral forces at the rear and the front wheel, modeled using the Pacejka tire model \cite{MF}. For more details and the exact tire model formulation, we refer to \cite{Liniger_2014}.

Based on the bicycle model \eqref{eq:odeModel}, the path planning model can be derived by looking at constant velocity points. These are points in the model where $\dot{v}_x = \dot{v}_y = \dot{\omega} = 0$. By selecting $N_m$ such points, in our case $N_m = 129$ a library of motion primitives can be formed. If the velocity is fixed to one of these constant velocities the three velocity states can be neglected, and the resulting trajectory is either a circle or a straight.

Based on these motion primitives one can build a new vehicle model, where the car drives with one of these constant velocities for a fixed time $T_{pp}$, after that time a new constant velocity point can be selected from the library. This direct concatenation of these constant velocity segments is motivated by the time scale separation present in the system where the velocities change fast enough, such that at a sampling time of $T_{pp}$ they can be regarded as constant. However, for this assumption to hold we need to constrain the motion primitives that can be concatenated. We do this by adding an automaton to the dynamics which limits the concatenation of constant velocity points. The automaton is generated such that there exist inputs to the bicycle model \eqref{eq:odeModel} which achieve the transition from one to another constant velocity point within a time that is smaller than $T_{pp}$. For our model, we use a sampling time of $T_{pp} = 0.16$\,s, whereas the transition can be achieved within 0.1\,s. Note that the model is hybrid due to the discrete mode, which corresponds to the constant velocity point, and the continuous dynamics describing the position and heading when the car drives with a specific constant velocity mode. However, when the continuous trajectory between the sampling points, where the constant velocity point can be switched, is neglected. The system can be formulated as a discrete-time system by embedding the discrete mode $q$ in the real line. The resulting discrete time system is given by $\bar{x}_{k+1} = \bar{f}(\bar{x},u_k)$ with $u \in U(q)$, as described in Section \ref{sec:GameIngreDyn}. Where $\bar{f}(\bar{x},u_k)$ describes the kinematic part of the car model when driving with a constant velocity (integration of the position and heading dynamics in \eqref{eq:odeModel}) as well as the dynamics on the discrete mode $q$, and $u \in U(q)$ is the automaton concatenation constraint.
 
 \end{appendices}

\section*{Acknowledgment}
The authors would like to thank Damian Frick for helpful discussions and advice.

\bibliographystyle{IEEEtran}
\bibliography{bimat_racing}
    
\end{document}